\documentclass[12pt]{amsart}
\usepackage{amscd}
\usepackage{amsfonts}
\usepackage{bbm}
\usepackage{times}
\usepackage{amssymb}
\usepackage{mathrsfs}
\usepackage{tabularx}
\usepackage{amsthm}
\usepackage{amsmath}
\usepackage{hyperref}

\numberwithin{equation}{section}

\begin{document}

\title[On the Lipschitz continuity of certain quasiregular mappings]
{On the Lipschitz continuity of certain quasiregular mappings between smooth Jordan domains}
\author{Jiaolong Chen, Peijin Li, Swadesh Kumar Sahoo and Xiantao Wang*}

\address{
Department of Mathematics, Shantou University, Shantou, Guangdong 515063, People's Republic of China}
\email{
jiaolongchen@sina.com (J. Chen)}

\address{
Department of Mathematics, Indian Institute of Technology Indore, Indore 452 017, India}
\email{swadesh@iiti.ac.in (S. Sahoo)}

\address{Department of Mathematics, Hunan Normal University, Changsha, Hunan 410081, People's Republic of China}
\email{wokeyi99@163.com (P. Li)}

\address{
Department of Mathematics, Shantou University, Shantou, Guangdong 515063, People's Republic of China}
\email{xtwang@stu.edu.cn (X. Wang)}

\thanks{The research was partly supported by NSF of China (No. 11571216)} 
\begin{abstract} We first investigate the Lipschitz continuity of $(K, K')$-quasiregular $C^2$ mappings between two Jordan domains with
smooth boundaries, satisfying certain partial differential inequalities concerning Laplacian.  Then two applications of the obtained result are given:
As a direct consequence, we get the Lipschitz continuity of $\rho$-harmonic $(K, K')$-quasiregular mappings, and as the other application, we study the Lipschitz continuity of $(K,K')$-quasiconformal self-mappings of the unit disk, which
are the solutions of the Poisson equation $\Delta w=g$. These results generalize
and extend several recently obtained results by Kalaj, Mateljevi\'{c} and Pavlovi\'{c}.
\thanks{*Corresponding author.}

\vskip 3mm
\end{abstract}

\maketitle

%

\section{Introduction and statement of the main results}\label{csw-sec1}
\subsection{Preliminaries}
\subsubsection{$(K, K')$-quasiregular mappings}

Let $$
A= \left(
   \begin{array}{cc}
     a & b \\
     c & d \\
   \end{array}
   \right)\in \mathbb{R}^{2\times2}.
 $$
We will consider the matrix norm:
$$|A|=\max\{|Az|:\; z\in \mathbb{C}, |z|=1\}$$
and the matrix function:
$$l(A)=\min\{|Az|: \;z\in \mathbb{C}, |z|=1\}.$$
Let $D$ and $\Omega$ be subdomains of the complex plane $\mathbb{C}$, and let $w=u+iv:$ $D\to \Omega$ be a function
that has both partial derivatives at $z=x+iy$ in $D$. $\nabla w$ denotes the Jacobian matrix
$$
\left(
  \begin{array}{cc}
    u_x & u_y \\
    v_x & v_y \\
  \end{array}
\right).$$
Obviously,
$$|\nabla w|=|w_z|+|w_{\overline z}|\;\;\mbox{and}\;\; l(\nabla w)=\big ||w_z|-|w_{\overline z}|\big |.$$

We say that a function $w: D\to \mathbb{C}$ is {\it absolutely continuous on lines}, {\it ACL} in brief, in the region $D$ if for
every closed rectangle $R\subset D$ with sides parallel to the axes $x$ and $y$, $w$ is absolutely continuous on almost every horizontal line
and almost every vertical line in $R$. Such a function has, of course, partial derivatives $w_x$ and $w_y$ a.e. in $D$.
Further, we say $w\in ACL^{2}$ if $w\in ACL$ and its partial derivatives are locally $L^2$ integrable in $D$.

We say that a sense-preserving continuous mapping $w:$ $D\to \Omega$ is
{\it  $(K, K')$-quasiregular} if

 (1) $w$ is $ACL^{2}$ in $D$ and $J_{w}\not=0$ a.e. in $D$;

(2) there are constants $K\geq1$ and $K'\geq0$ such that
 $$|\nabla w|^2\leq KJ_w+K',
 $$
where $J_w$ denotes the Jacobian of $w$, which is given by
$$J_w=|w_z|^2-|w_{\overline z}|^2=|\nabla w|l(\nabla w).$$

In particular, if $w$ is a $(K, K')$-quasiregular homeomorphism,
then $w$ is $(K, K')$-quasiconformal;
and if $w$ is a $K$-quasiregular homeomorphism, then $w$ is $K$-quasiconformal.

If $K'=0$, then ``$(K, K')$-quasiregular'' (resp. ``$(K, K')$-quasiconformal") mappings reduce to ``$K$-quasiregular'' (resp. ``$K$-quasiconformal").

We remark that there are $(K, K')$-quasiregular mappings which are not $K_1$-quasiregular for any $K_1\geq 1$, and also there are
$(K, K')$-quasiconformal mappings whose inverses are not $(K_1, K_1')$-quasiconformal for any $K_1\geq 1$ and $K_1'\geq 0$. See the examples in Sections \ref{sec-2} and \ref{sec-4} for the details.

The following result easily follows from the proof of Lemma $2.1$ in \cite{chen}.

\newtheorem*{lem1.1}{Lemma 1.1}
\begin{lem1.1}
Suppose $w$ is a $(K, K')$-quasiregular mapping.
Then,
$$|\nabla w| \leq K l(\nabla w)+\sqrt{K'}.$$
\end{lem1.1}

A mapping $f:$ $D\to \Omega$ is {\it proper} if the preimage of every compact set in $\Omega$ is compact in $D$. It is known that if $D=\Omega=\mathbb{D}=\{z:|z|<1\}$, then the mapping $f$ is proper if and only if $|f(z)|\rightarrow 1$ as $|z|\rightarrow 1$ (cf. \cite{KP}).

\subsubsection{Lipschitz continuity}

We say that a mapping $f:$ $D\to \Omega$ is in $Lip_{\alpha}$
if there exists a constant $L_{1}$ and an exponent $ \alpha\in (0,1] $ such that for all $z,$ $w\in D$,
$$|f(z)-f(w)|\leq L_{1}|z-w|^{\alpha}.$$
Such mappings are also called {\it $\alpha$-H\"{o}lder continuous}.

In particular, if $\alpha=1$, then we say that $f$ is {\it Lipschitz continuous}.

The mapping $f$ is said to be {\it coLipschitz continuous} if there exists a constant $L_{2}$ such that for all $z,$ $w\in D$,
$$|f(z)-f(w)|\geq L_{2}|z-w| .$$

\subsubsection{Jordan domains}
A {\it Jordan curve} is a set in the complex plane $\mathbb{C}$ which is homeomorphic to a circle. For a {\it Jordan domain}, we mean a domain whose boundary is a Jordan curve. In this paper, unless specially stated, all Jordan domains are assumed to be bounded.

Denote by $\ell(\gamma)$ the length of $\gamma$, and let $\Gamma: [0, \ell(\gamma)]\rightarrow\gamma$ be the arc length parameterization of $\gamma$, i.e. the parameterization satisfying the condition:
\begin{center}
 $|\Gamma'(s)|=1$ \;\;a.e. in $[0,\ell(\gamma)]$.
\end{center}

We say that $\gamma$ is of {\it class $C^{n,\alpha}$} for some $n\in \mathbb{N}$
and $ \alpha\in(0,  1]$ if $\Gamma$ is of class $C^n$ and
$$\sup_{t\not=s\in[0,\ell(\gamma)]} \frac{|\Gamma^{(n)}(t)- \Gamma^{(n)}(s)|}{|t-s|^{\alpha}}<\infty,$$
where $\Gamma^{n}(t)$ denotes the $n^{th}$ derivative of $\Gamma (t)$ with respect to $t$.
The Jordan domain $D$ is called a {\it $C^{n,\alpha}$ domain} if its
 boundary $\partial D$ is a $C^{n,\alpha}$ Jordan curve.

Let $\gamma \in C^{1,\alpha}$ be a closed Jordan curve, and $d_{\gamma}(\Gamma(s),\Gamma(t))$ the distance
between $\Gamma(s)$ and $\Gamma(t)$ along the curve $\gamma$, i.e.
$$d_{\gamma}(\Gamma(s),\Gamma(t))=\min \{|s-t|,\ell(\gamma)-|s-t|\}.$$

A closed rectifiable Jordan curve $\gamma $ is said to enjoy a {\it $b$-chord-arc
condition} if there exists $b>1$ such that for all $z_{1}$, $z_2 \in \gamma$,
$$d_{\gamma}(z_1,z_2)\leq b|z_1 - z_2|.$$

We remark that the unit circle $\mathbb{S}$ enjoys a $\frac{\pi}{2}$-chord-arc condition.
\medskip

\subsubsection{Normalized mappings}
For a closed curve $\beta$ in $\mathbb{C}$, three points $a_0$, $a_1$ and $a_2$ in $\beta$ are said to be {\it well-distributed} if for $i\in \{0,1\}$,
$$\ell(\beta[a_i,a_{i+1}])=\ell(\beta[a_{i+1},a_{i+2}]),$$ where $a_3=a_0$ and $\beta[a_i,a_{i+1}]$ denotes the part of $\beta$ with endpoints $a_i$ and $a_{i+1}$.

Let $D$ be a Jordan domain with rectifiable boundary. We will say that a
mapping $f:$ $\overline{\mathbb{D}}=\mathbb{D}\cup \mathbb{S}\to \overline{\Omega}$ is {\it normalized}
if there are three well-distributed points $t_0$, $t_1$, $t_2$ in $\mathbb{S}$, their images $f(t_0)$, $f(t_1)$ and $f(t_2)$ under $f$ are also well-distributed in $\partial \Omega=f(\mathbb{S})$ (cf. \cite{KM}).

\subsection{Lipschitz continuity for certain $(K, K')$-quasiregular mappings}

If a $(K, K')$-quasiregular (resp. $(K, K')$-quasiconformal) mapping is harmonic, then it is said to be {\it harmonic $(K, K')$-quasiregular} (resp. {\it harmonic $(K, K')$-quasiconformal}).  Martio \cite{M} was the first who considered harmonic quasiconformal mappings
in $\mathbb{C}$.
Recent papers \cite{chen, K1, KM0, KM, KP0, LCW, P} and references therein together bring much light
on the topic of harmonic quasiconformal mappings in $\mathbb{C}$. See \cite{K2, Kalaj}
for the discussions in this line in the space. In \cite{MKM, MV},
the Lipschitz characteristic of harmonic quasiconformal mappings  has been discussed. See \cite{TW1, TW2, TW3, W} for similar discussions in this line.
In \cite{roja}, Finn and Serrin discussed the H\"{o}lder continuity of a class of elliptic mappings which satisfy the following partial differential inequality: $$|w_{z}|^{2}+|w_{\overline{z}}|^{2} \leq K J_{w} +K',$$
where $K\geq 1$ and $K'\geq 0$ are constants. See also \cite{Nir}.

Recently, Kalaj and Mateljevi\'{c} \cite{KM} discussed the Lipschitz continuity of $(K, K')$-quasiconformal harmonic mappings. In \cite{KM0}, they
considered the Lipschitz continuity of a quasiconformal $C^2$
diffeomorphism $w:$ $D\to \Omega$ which satisfies the partial differential inequality:
$$|\Delta w|\leq M|w_z\cdot w_{\overline{z}}|,\eqno{(1.1)}$$
where $M\geq0$ is a constant and $D$ (resp. $\Omega$) denotes a $C^{1,\alpha}$ (resp. $C^{2,\alpha}$) Jordan domain,
and in \cite{K3}, as a generalization of the discussions in \cite{KM0},
Kalaj studied the Lipschitz continuity of a $K$-quasiregular $C^2$
mapping $w:$ $D\to \Omega$ which satisfies the partial differential inequality:
$$|\Delta w|\leq M|\nabla w|^2+N \eqno{(1.2)}$$
for some constants $M\geq 0$ and $N\geq 0$.

Obviously, if  a mapping satisfies the partial differential inequality (1.1), then it also satisfies (1.2).
Observe that if $M=N=0$ in (1.2), then $w$ is harmonic. The reader is referred to \cite{D} for the properties of this class of mappings.

As the first aim of this paper, we consider the Lipschitz continuity of $(K, K')$-quasiregular $C^2$ solutions of (1.2).
Our result is as follows.
\newtheorem*{thm1.1}{Theorem 1.1}
\begin{thm1.1} Suppose $w$ is a proper $(K, K')$-quasiregular $C^2$ mapping of a Jordan domain $D$
with $C^{1,\alpha}$ boundary onto a Jordan domain $\Omega$ with $C^{2,\alpha} $ boundary. If $w$ satisfies the partial differential inequality (1.2) for constants $M>0$ and $N\geq 0$,
then $w$ has bounded partial derivatives in $D$.
In particular, $w$ is Lipschitz continuous.
\end{thm1.1}

We remark that Theorem 1.1 is a substantial
generalization of \cite[Theorem 1.2]{K3} and \cite[Theorem 1.3]{KM0} (The example in Section \ref{sec-2} demonstrates this). This theorem will be proved in Section \ref{sec-3}.

\subsection{Lipschitz continuity for certain $\rho$-harmonic $(K, K')$-quasiregular mappings}\label{ssec-4}

A two times continuously differentiable complex-valued function $f$ between two domains $D$ and $\Omega$ in $\mathbb{C}$ is said to be  {\it $\rho$-harmonic} if it satisfies the Euler-Lagrange equation
$$f_{z\overline{z}}+\big((\log \rho)_{w}\circ w\big) f_{z}f_{\overline{z}}=0,\eqno{(1.3)}$$
where $w=f(z)$, and $\rho(w)|dw|$ is an arbitrary conformal $C^1$-metric defined in $\Omega$.

 If $\varphi$ is a holomorphic mapping different from 0 and if $\rho=|\varphi|$ in $\Omega$, we
call $w$ a $\varphi$-harmonic mapping.

Since $\rho^2=\varphi \overline{\varphi}$, an elementary computation yields $2(\log \rho)_w=(\log \varphi)'$.
It follows from (1.3) that if $f$ is $\varphi$-harmonic, then
$$f_{z \overline{z}}+\left(\frac{\varphi'}{2 \varphi}\circ w\right) f_{z}f_{\overline{z}}=0.\eqno{(1.4)}$$

In \cite{KM0}, Kalaj and Mateljevi\'{c} considered the Lipschitz continuity of $\rho$-harmonic quasiconformal mappings.
As a direct consequence of Theorem 1.1, we obtain the Lipschitz continuity of $\rho$-harmonic $(K, K')$-quasiregular mappings.
Our first result concerns the Lipschitz continuity of $\varphi$-harmonic $(K,K')$-quasiregular mappings, which is as follows.
\newtheorem*{cor1.1}{Corollary 1.1}
\begin{cor1.1}
Suppose $w$ is a $\varphi$-harmonic mapping of $\mathbb{D}$ onto a $C^{2,\alpha}$
Jordan domain $\Omega$. If $|(\log \varphi)'|_{\infty}=\sup\{|(\log \varphi(z))'|:z\in \mathbb{D}\} <\infty$ and $w$ is a proper $(K,K')$-quasiregular mapping, then $w$ has
bounded partial derivatives. In particular, $w$ is Lipschitz continuous.
\end{cor1.1}

Our next result, concerning approximately analytic metrics, generalizes Corollary 1.1, where a $C^1$ function $h$ is said to be {\it approximately analytic} if there is a constant $C\geq 0$ such that $|h_{\overline{z}}|\leq C|h|$.

\newtheorem*{cor1.2}{Corollary 1.2.}
\begin{cor1.2} Suppose $w$ is a $\rho$-harmonic mapping of $\mathbb{D}$ onto a $C^{2,\alpha}$
Jordan domain $\Omega$. Further, if $\rho$ is approximately analytic in $\Omega$ and $w$ is a proper $(K,K')$-quasiregular mapping, then $w$ has
bounded partial derivatives. In particular, $w$ is Lipschitz continuous.
\end{cor1.2}

We remark that Corollary 1.1 (resp. Corollary 1.2) is a
generalization of \cite[Theorem 3.1]{KM0} (resp. \cite[Theorem 3.3]{KM0}).

\subsection{Lipschitz continuity for a class of $(K, K')$-quasiconformal self-mappings of the unit disk}\label{sec-5}

In this subsection, we discuss the Lipschitz continuity of $(K, K')$-quasiconformal solutions of the Poisson equation (1.5) (see below) in $\mathbb{D}$. We start with some necessary definitions.
Let $P$ be the Poisson kernel, i.e. the function
$$P(z,e^{i\theta})=\frac{1-|z|^2}{|z-e^{i\theta}|^2},$$
and let $G$ denote the Green function of $\mathbb{D}$, i.e.
$$G(z,\omega)=\frac{1}{2\pi}\log\Big|\frac{1-z\overline\omega}{z-\omega}\Big|,$$
where $z\in \mathbb{D}\setminus\{\omega\}$. Obviously, $P$ is harmonic in $\mathbb{D}$ (cf. \cite{D}), and
$G$ is harmonic in $\mathbb{D}\setminus\{\omega\}$.

Let $f:$ $\mathbb{S}\rightarrow\mathbb{C}$ be a bounded integrable function in $\mathbb{S}$,
and let $g:$ $\mathbb{D}\rightarrow\mathbb{C}$ be continuous. It is known that the solutions of the Poisson equation
$$\Delta w=g\eqno{(1.5)}$$
in $\mathbb{D}$ satisfying the boundary condition $w|_{\mathbb{S}}=f\in L^1(\mathbb{S})$ has the following expression:
$$w =P[f]-G[g],$$
where
$$P[f](z)= \frac{1}{2\pi}\int_0^{2\pi}P(z,e^{i\varphi})f(e^{i\varphi})d\varphi,$$
$$G[g](z)=\int_{\mathbb{D}}G(z,\omega)g(\omega){\rm dm}(\omega)$$
and ${\rm dm}(\omega)$ denotes the Lebesgue measure in $\mathbb{C}$. Also, it is known that if $f$ and $g$ are continuous
in $\mathbb{S}$ and $\overline{\mathbb{D}}$, respectively, then $w$ has a continuous extension $\tilde{w}$ to $\mathbb{S}$ and $\tilde{w}|_{\mathbb{S}}=f$ (cf. \cite{LH}).

For convenience, in the following, we always set
$$\mathcal{P}=P[f]\;\;\mbox{and}\;\;\mathcal{G}=G[g].$$

Let $\mathcal{D}_{\mathbb{D}}(g)$ denote the family of all solutions $w$ of the Poisson equation (1.5) from $\mathbb{D}$ onto $\mathbb{D}$, which satisfy
that each element $w$ is a $C^2$ diffeomorphism,
each restriction $w|_{\mathbb{S}}=f$ is normalized, every function $f(e^{it})=e^{i \psi(t)}$
is an absolutely homeomorphism of $\mathbb{S}$ onto $\mathbb{S}$ and $\psi(2\pi)=\psi(0)+2\pi$.

Observe that any element in $\mathcal{D}_{\mathbb{D}}(g)$ is proper and satisfies the partial differential inequality (1.2) with $M=0$ and $N=|g|_{\infty}$.

In \cite{KP}, Kalaj and Pavlovi\'{c} discussed the Lipschitz continuity of quasiconformal self-mappings of $\mathbb{D}$ satisfying the Poisson equation (1.5).
As the main application of Theorem 1.1, we study the Lipschitz continuity of $(K, K')$-quasiconformal solutions of (1.5).
The aim is to generalize the arguments in \cite{KP} to the case of $(K,K')$-quasiconformal solutions of (1.5).
 The following is our result.

\newtheorem*{thm1.2}{Theorem 1.2.}
\begin{thm1.2} Suppose that $w \in \mathcal{D}_{\mathbb{D}}(g)$ is a $(K,K')$-quasiconformal mapping and that $g\in C(\overline{\mathbb{D}})$.

 \noindent $(1)$
 Then for all $z_1$ and $z_2\in \mathbb{D}$, $$|w(z_1)-w(z_2)|\leq M|z_1 -z_2|,$$
 where $M=M(K, K', |g|_{\infty})$ which means that the constant $M$ depends only on $K$, $K'$ and $|g|_{\infty}$.
 \medskip

 \noindent $(2)$ If $w^{-1}$ is also $(K,K')$-quasiconformal, then for all $z_1$ and $z_2\in \mathbb{D}$, $$|w(z_1)-w(z_2)|\geq N|z_1 -z_2|,$$
 where $N=N(K, K', |g|_{\infty})$.
\end{thm1.2}

By comparing with \cite[Theorem 1.2]{KP}, a natural question is that whether the assumption ``$w^{-1}$ being $(K,K')$-quasiconformal" in the second statement in Theorem
 1.2 is necessary or not. We will construct an example (Example 2.1 below) to show that there is a solution of the Poisson equation such that it is $(K, K')$-quasiconformal, its inverse is not $(K,K')$-quasiconformal for any $K\geq 1$ and $K'\geq 0$, and it is not coLipschitz continuous. This fact shows that the mentioned assumption in Theorem 1.2 is necessary. Section 4 is devoted to the proof of Theorem 1.2 together with the statement and the proof of Example 4.1.

In Section \ref{sec-2}, we will construct an example to show that Theorem 1.1 is a substantial generalization of the corresponding results in \cite{K3} and \cite{KM0}, respectively.

\section{An example}\label{sec-2}
In this section, we will construct an example to show the existence of the $(K, K')$-quasiregular solutions of the partial differential inequality (1.2), which satisfy the requirements in Theorem 1.1, but fail to satisfy the assumptions in the corresponding results in \cite{K3} and \cite{KM0}, respectively.

\newtheorem*{exa2.1}{Example 2.1}
\begin{exa2.1}
Let $w(z)=2|z|^4z^2-|z|^{10}z^2$ in $\mathbb{D}$. Then
\begin{enumerate}
\item
$w$ is a $(1, 144)$-quasiregular mapping of $\mathbb{D}$
onto $\mathbb{D}$;
\item
$w$ is proper;
\item
$w$ is not $K$-quasiregular for any $K\geq1$;
\item
$w$ satisfies $|\Delta w|\leq |\nabla w|^2+76$;
\item
$w$ doesn't satisfy $|\Delta w|\leq M_1|\nabla w|^2$ for any $M_1\geq 0$;
\item
$w$ doesn't satisfy $|\Delta w|\leq M_2|w_z\cdot w_{\overline{z}}|$ for any $M_2\geq 0$;
\item
$w$ is Lipschitz continuous.
\end{enumerate}
\end{exa2.1}
\begin{proof} Obviously, $w(e^{i\theta})=e^{2i\theta}$ for $\theta\in [0, 2\pi]$,
 $$w_z= 8|z|^4 z-7|z|^{10} z\;\; \mbox{and}\;\; w_{\overline{z}}=4|z|^2 z^3-5|z|^8 z^3.$$
 It follows that
 $$
|w_z(z)|-|w_{\overline z}(z)|=\begin{cases}
\displaystyle \,2|z|^5(2- |z|^{6})>0,\;\;\;\text{if}\;\; 0<|z|^6 <4/5,\\
\displaystyle12|z|^5(1-|z|^{6})>0,\;\;\text{if}\;\; 4/5\leq|z|^6<1
\end{cases}$$
in $\mathbb{D}\setminus \{0\}$, and
$$
|\nabla w |^2-J_{w} \leq |\nabla w |^2=\begin{cases}
\displaystyle 144|z|^{10}(1-|z|^{6})^2<144,\;\;\text{if}\;\; 0<|z|^6 <4/5,\\
\displaystyle \;\;\; 4|z|^{10}(2-|z|^{6})^2<16,\;\;\;\;\,\text{if}\;\;  4/5 \leq |z|^6<1.
\end{cases} $$
Let $z_1=re^{i \theta}$ and $z_2=re^{i (\theta+\pi)}$ with $0<r<1$. Then $z_1\not=z_2\in \mathbb{D}$ and $w(z_1)=w(z_2)$. Hence we have proved that $w$ satisfies the first two assertions in the example.

The limits
$$\lim_{|z|\to1^-}\frac{|w_{\overline z}|}{|w_z|}=\lim_{|z|\to1^{-}}\frac{ 5|z|^{11}  - 4|z|^5}{8|z|^5-7|z|^{11}}=1$$
and
$$\lim_{|z|\to1^-}(|w_z(z)|+|w_{\overline z}(z)|)=2$$
tells us that $w$ is not $K$-quasiregular for any $K\geq 1$, which implies that the third assertion is satisfied.

Obviously, $$|\Delta w|=\big|64|z|^4-140|z|^{10}\big|\leq 76\leq |\nabla w|^2+76$$
and
$$\lim_{|z|\to0}\frac{|\Delta w|}{|\nabla w|^2}=\lim_{|z|\to0}\frac{16-35|z|^{6}}{36|z|^{6}(1- |z|^{6})^2}=+\infty.$$
Hence the assertions from forth to sixth in the example hold.

It remains to show that the Lipschitz continuity of $w$. This easily follows from the estimate $|\nabla w|<12$ in $\mathbb{D}$, and so the proof of the example is complete.
\end{proof}

\section{Lipschitz continuity for certain $(K, K')$-quasiregular mappings between Jordan domains}\label{sec-3}

The main purpose of this section is to prove Theorem 1.1. Before the proof, some preparation is needed.

\subsection{Some auxiliary results} We start this subsection with a lemma.

\newtheorem*{lem3.1}{Lemma 3.1}
\begin{lem3.1}
Suppose that $w$ is a $(K, K')$-quasiregular mapping from $D$ to $\Omega$, where both $D$ and $\Omega$ are Jordan domains in $\mathbb{C}$, and that $w$ has the decomposition $w=\rho S$, where $\rho=|w|$. Then we have
$$|\nabla w|\leq K|\nabla \rho|+\sqrt{K'}\eqno{(3.1)}$$
and
$$
\frac{|\nabla \rho|-\sqrt{K'}}{K}\leq \rho|\nabla S|\leq K|\nabla \rho|+\sqrt{K'}\eqno{(3.2)}$$
a.e. in $D$.
\end{lem3.1}
\begin{proof}
 First, we prove the inequality (3.1). Since $\rho=|w|$, similar discussions as in the proof of \cite[Lemma 2.2]{K3} guarantee that
$$\nabla \rho =\frac{1}{|w|}w^{T}\nabla w,$$
where $T$ denotes the transpose of matrices. Here and hereafter, we regard $w=u+iv$ as not only a number in $\mathbb{C}$, but also a $2\times 1$
column vector, where both $u$ and $v$ are real. Then
$$\nabla \rho h=\frac{1}{|w|}\langle\nabla w h,w\rangle, \eqno{(3.3)}$$ and so
$$ |\nabla \rho|=\max_{|h|=1}|\nabla \rho h|=\max_{|h|=1}\frac{|\langle\nabla w h,w\rangle|}{|w|}\leq |\nabla w|. \eqno{(3.4)}$$
It follows from Lemma 1.1 that to prove (3.1), it suffices to show $$l(\nabla w)\leq |\nabla \rho|.$$

Since $J_{w}=\det(\nabla w)\not=0$ a.e. in $D$, we know that there exists $h_1$ such that
$$\nabla w h_1=\frac{w}{|w|} \;\,\mbox{a.e. in}\;\, D.\eqno{(3.5)}$$
Then we infer from (3.3) that $$\nabla \rho h_1= \frac{< \frac{w}{|w|},w>}{|w|}=1,$$ and so for $h= \frac{h_1}{|h_1|}$, it easily follows that $\nabla \rho h= \frac{1}{|h_1|}$.
Also by (3.5), we get $h_1=(\nabla w)^{-1} \frac{w}{|w|}$ a.e. in $D$, and thus
$$|\nabla \rho|\geq |\nabla \rho h|=\frac{1}{|h_1|}=\left|(\nabla w)^{-1} \frac{w}{|w|}\right|^{-1}.$$
Since $$\left|(\nabla w)^{-1} \frac{w}{|w|}\right|\leq |(\nabla w)^{-1} | $$
and
$$|(\nabla w)^{-1}|=\frac{1}{\min\{| \nabla w h|:|h|=1\}},$$
we obtain that
$$|\nabla \rho|\geq \left|(\nabla w)^{-1} \right|^{-1}=\min\{|\nabla w h|:|h|=1\}=l(\nabla w),\eqno{(3.6)}$$
as required. Hence the inequality (3.1) holds.\medskip

Now, we check the truth of the right side of the inequality (3.2).
We infer from a similar argument as in the proof of \cite[Lemma 2.2]{K3} that
$$\nabla S=\frac{\nabla w}{|w|}-\big((\nabla w)^{T} w\big )\otimes \frac{w}{|w|^3},$$
and so for all $h\in \mathbb{C}$,
$$ \nabla S h=\frac{\nabla w h}{|w|}- \frac{w\langle\nabla w h, w \rangle}{|w|^3},$$
  where $\otimes$ is the tensor product between column vectors, i.e. for two vectors $\overrightarrow{a}$ and $\overrightarrow{b}$, $\overrightarrow{a}\otimes\overrightarrow{b}=\overrightarrow{b}\cdot (\overrightarrow{a})^{T}$.
Then
$$\rho^2 |\nabla Sh|^2=|\nabla wh|^2- \Big\langle\nabla wh,\frac{w}{|w|} \Big\rangle^2 \leq |\nabla w|^2|h|^2,
\eqno{(3.7)}$$
whence
$$\rho |\nabla S |=\sup_{|h|\not=0}\left\{\rho\left|\nabla S \frac{h}{|h|}\right|\right\} \leq |\nabla w|.$$
Hence we obtain from Lemma 1.1 and (3.6) that
$$\rho|\nabla S|\leq |\nabla w|\leq K l(\nabla w)+\sqrt{K'}\leq K|\nabla \rho|+\sqrt{K'}, $$ as needed.
\medskip

Next, we check the truth of the left side of (3.2).
Since, obviously, there is an $h_1$ such that $\langle\nabla w h_1 ,\frac{w}{|w|}\rangle=0$ with $|h_1|=1 $, we see from (3.7) that
$$\rho^2|\nabla S h_{1}|^2=|\nabla w h_1|^2,$$
which, together with Lemma 1.1, guarantees that
$$\rho |\nabla S|\geq \rho \left|\nabla S \frac{h_1}{|h_1|}\right|=\left|\nabla w \frac{h_1}{|h_1|}\right|\geq l(\nabla w)\geq\frac{|\nabla w|-\sqrt{K'}}{K},$$
and then we get from (3.4) that
$$\rho |\nabla S|\geq \frac{|\nabla w|-\sqrt{K'}}{K}\geq \frac{|\nabla \rho|-\sqrt{K'}}{K},$$ which is what we want. Therefore, the proof of the lemma is complete.
\end{proof}

\vspace{0.3truecm}
\noindent{\bf Remark 3.1.} Lemma 3.1 is a generalization of \cite[Lemma 2.2]{K3}.

The proof of Theorem 1.1 also needs the following known results.

\newtheorem*{lema}{Lemma A (\cite[Proposition 2.5]{K3} or \cite{go, War1, War2, Pom})}
\begin{lema}  If $D$ and $\Omega$ are Jordan domains having $C^{n,\alpha}$ $(n\geq1)$ boundary and if $\omega$ is a conformal mapping of $D$ onto $\Omega$, then

(a) $|\omega'(z)|\geq \inf\{|\omega'(\zeta)|:\zeta\in D\}>0$ for $z\in D$.

(b) $\omega^{(n)}\in C^{\alpha}(\overline{D})$. In particular $|\omega^{(n)}|_{\infty}=\sup\{|\omega^{n}(z)|:z\in D\}<\infty$.
\end{lema}

\newtheorem*{lemb}{Lemma B (\cite[Theorem 4']{Hei} or \cite[Proposition 3.2]{KM1})}
\begin{lemb}
Suppose $u$ is a continuous function from $\overline{\mathbb{D}}$ into the real axis $\mathbb{R}$ and satisfies the following conditions:

(1) $u$ is $C^2$ in $\mathbb{D}$;

(2) $u_b(\theta)=u(e^{i\theta})$ is $C^2$; and

(3) $|\Delta u|\leq M_0 |\nabla u|^2+N_0$ in $\mathbb{D}$ for some constants $M_0$ and $N_0$.

Then $|\nabla u|$ is bounded in $\mathbb{D}$.
\end{lemb}

\subsection{The proof of Theorem 1.1}

 We are ready to prove the Lipschitz continuity of $w$. We divide the proof into two steps. In the first step, we construct a $(K_1, K'_1)$-quasiregular self-mapping $ \tau$ in $\mathbb{D}$ satisfying the partial differential inequality (1.2) for some constants $M_1$ and $N_1.$ In the second step, by applying the mapping $\tau$, we show that $|\nabla w|$ is bounded in $D$, which completes the proof.

\vspace{0.3truecm}
\noindent{\bf Step 3.1.} The construction of $\tau$.

Let $\varphi$ be a conformal mapping of $\mathbb{D}$ onto $D$, and $g$ a conformal mapping of $\Omega$ onto $\mathbb{D}$, and $$\tau=g\circ w \circ\varphi.$$

In the following, we apply Lemma A to show that this $\tau$ satisfies our requirements as mentioned in the first paragraph in this subsection.

First, it follows from
$$|\nabla\tau|=|\tau_z|+|\tau_{\overline{z}}|\;\;\mbox{and}\;\;|\nabla\tau|=|g'\varphi'||\nabla w|\eqno{(3.8)}$$  that
$$J_{\tau}=|g'\varphi'|^2J_w.$$
Meanwhile, by
Lemma A, we know that
the function $|g'|$ (resp. $|\varphi'|$) is bounded from above and below. Hence
we easily know that $\tau$ is a $(K, K'|g'|^2_{\infty} |\varphi'|_{\infty}^2)$-quasiregular self-mapping in $ \mathbb{D}$.

Second, since
$$|4\tau_z \tau_{\overline{z}}|\leq(|\tau_z|+|\tau_{\overline{z}}|)^2=|\nabla \tau|^2,$$
$$ |\Delta w|\leq M|\nabla w|^2+N=M\frac{|\nabla \tau|^2}{|g'|^2|\varphi'|^2}+N,$$ and since elementary computations lead to
$$\Delta\tau=(4g'' w_{z}w_{\overline{z}}+ g' \Delta w)|\varphi'|^2=\Big(\frac{4g'' \tau_{z}\tau_{\overline{z}}}{g'^2}+g'|\varphi'|^2\Delta w\Big),$$
 we have
$$ |\Delta \tau|\leq \Big(\frac{|g''|_{\infty}}{|g'|_1^2}+ \frac{M}{|g'|_1}\Big)|\nabla \tau|^2+|g'|_{\infty} |\varphi'|^2_{\infty}N,\eqno{(3.9)}$$
where $|g'|_1=\inf\{|g'(z)|:\; z\in G\}$. The boundedness of $|g''|_{\infty}$ follows from Lemma A. Hence $\tau$ satisfies (1.2), and so this $\tau$ is what we need.\medskip

\vspace{0.3truecm}
\noindent{\bf Step 3.2.}  $|\nabla w|$ is bounded in $D$.

It follows from (3.8) that $$|\nabla w|=|\nabla \tau| |g'\varphi'|^{-1},$$ and then we know from Lemma A that to prove the boundedness of $|\nabla w|$ in $D$, it suffices to show the boundedness of $|\nabla \tau|$ in $\mathbb{D}$. Now, we are going to prove the boundedness of $|\nabla \tau|$.
Obviously, it follows from the fact $w\in C^2(\mathbb{D})$ that $\tau\in C^2(\mathbb{D})$ and then $|\nabla \tau |$ is bounded in $\overline{\mathbb{D}}_{r}$ for any $r\in (0,1)$, where $\overline{\mathbb{D}}_{r}=\{z\in \mathbb{C}:\; |z|\leq r\}$. To prove the boundedness of $|\nabla \tau |$ in $\mathbb{D}$, it is enough to show that $|\nabla \tau |$ is bounded in $\mathbb{D}\setminus\overline{\mathbb{D}}_{r}$ for some $r\in (0,1)$. For this,
we let $$\tau=\rho S,$$ where $\rho=|\tau|$. Then the inequality (3.1) in Lemma 3.1 makes sure that if $|\nabla\rho|$ is bounded, then so is $|\nabla  \tau|$. To prove that $|\nabla\tau|$ is bounded in $\mathbb{D}$, it is sufficient to find an $r\in (0,1)$ such that $|\nabla\rho|$ is bounded in $\mathbb{D}\setminus\overline{\mathbb{D}}_{r}$.

In the following, we apply Lemma B to show the existence of the needed $r$. To reach this aim, we need the following existence of a function related to $\rho$.

\vspace{0.3truecm}
\noindent{\bf Claim 3.1.} There is a function $\rho_{2}$ in $\mathbb{D}$ such that

(1)$\rho_{2}$ satisfies all assumptions in Lemma B;
 and

(2) $\rho_{2}$ and $\rho$ coincide with each other in $\mathbb{D}\setminus\overline{\mathbb{D}}_{r}$ for some $r\in (0,1)$.

We will apply Whitney's theorem \cite[Theorem 1]{wh} to construct such a function $\rho_{2}$. Since $w$ is a proper ($K,K'$)-quasiregular mapping, it follows that $\tau$ is a proper self-mapping of $\mathbb{D}$. Thus
$$\lim_{|z|\to 1^{-}}\rho(z)=1.$$  Therefore, there exists an $r_{1}>0$ such that $r_{1}\leq |z|\leq 1$ implies $$\rho(z)\geq 1/2.$$ Let $0<r_{1}< r_2<1$ and $r_2>1/2$.
 Since $\rho\in C^2(A)$, where $A=\{z:r_{1}\leq|z|\leq r_2\}$, according to Whitney's theorem, there exists an extension $\rho_1$ of the restriction $\rho|_{A}$ such that
$\rho_{1}\in C^2(\mathbb{C})$. Let

$$
\rho_{2}=\begin{cases}
\displaystyle\rho,\;\;z\in \mathbb{D}\setminus\overline{\mathbb{D}}_{\frac{r_1+r_2}{2}}, \\
\displaystyle \rho_{1},\;z\in\overline{\mathbb{D}}_{\frac{r_1+r_2}{2}}.
\end{cases}$$
Obviously, $\rho_{2}$ satisfies the assumptions (1) and (2) in Lemma B since $\rho_{2}(e^{i\theta})=\lim_{r\to 1^{-}}\rho_{2}(r e^{i\theta}) =\lim_{r\to 1^{-}}\rho(r e^{i\theta})=1$. In order to show that this $\rho_2$ is our needed, it remains to check that  $\rho_2$ satisfies (1.2), i.e. the third assumption in Lemma B. We will apply Lemma 3.1 to reach this goal.

It follows from \cite[Lemma 2.4]{K3} that
$$|\Delta \rho|=\Big|\frac{1}{2}\langle\Delta \tau, S\rangle+2 \rho |\nabla S|^2 \Big|
\leq \frac{1}{2}|\Delta \tau|+2\rho |\nabla S|^2.$$
Then the inequality (3.2) in Lemma 3.1 along with (3.9) implies
 \begin{eqnarray*}
|\Delta \rho| &\leq & \frac{1}{2}|\Delta \tau|+\frac{2}{\rho}\big(K|\nabla \rho|+\sqrt{ K'}\big)^2\\
&\leq& \frac{1}{2}\Big(\frac{|g''|_{\infty}}{|g'|_1^2}+ \frac{M}{|g'|_1  } \Big) |\nabla \tau|^2+\frac{1}{2}|g'|_{\infty}|\varphi'|^2_{\infty}N+
\frac{4}{\rho}(K^2|\nabla \rho|^2+K'),
\end{eqnarray*}
and further, we get from the inequality (3.1) in Lemma 3.1 that
$$
|\Delta \rho|
\leq N_1(\rho)|\nabla \rho|^2+M_1(\rho),
$$
where $$N_1(\rho)=\Big(\frac{4}{\rho}+\frac{|g''|_{\infty}}{|g'|_1^2}+ \frac{M}{ |g'|_1 } \Big)K^2$$
and
$$M_1(\rho)=\Big(\frac{4}{\rho}+\frac{|g''|_{\infty}}{|g'|_1^2}+ \frac{M}{ |g'|_1   } \Big)K'+\frac{1}{2}|g'|_{\infty}|\varphi'|^2_{\infty}N.$$
Since for all $z$ in $\mathbb{D}\setminus\overline{\mathbb{D}}_{r_1}$, $\rho(z)\geq \frac{1}{2}$, we see that
$$|\Delta \rho_{2}|\leq  N_1\big(\frac{1}{2}\big)|\nabla \rho_{2}|^2+ M_1\big(\frac{1}{2}\big)$$ in $\mathbb{D}\setminus\overline{\mathbb{D}}_{r_1}$.

Let
$$M_0=\max\left\{|\Delta \rho_{2}(z)|:\;z\in\overline{\mathbb{D}}_{\frac{r_1+r_2}{2}}\right\}.$$ Since $\rho_{1}\in C^2(\mathbb{C})$, by the definition of $\rho_2$,
we see that $M_0<+\infty$, and so
for all $z$ in $\mathbb{D}$,
$$|\Delta \rho_{2}|\leq  N_1\big(\frac{1}{2}\big)|\nabla \rho_{2}|^2+ M_1\big(\frac{1}{2}\big)+M_0,$$ which shows that $\rho_2$ satisfies the third assumption in Lemma B.

Since the definition of $\rho_2$ implies that $\rho_{2}=\rho$ in $\mathbb{D}\setminus\overline{\mathbb{D}}_{\frac{r_1+r_2}{2}}$, we see that this $\rho_2$ justifies our need and the proof of our claim is complete.
\medskip

It follows from Claim 3.1 and Lemma B that $|\nabla \rho_{2}|$ is bounded in $\mathbb{D}$. Hence $|\nabla \rho|$ is bounded in $\mathbb{D}\setminus\overline{\mathbb{D}}_{\frac{r_1+r_2}{2}}$, and so this radius $\frac{r_1+r_2}{2}$ is what is we wanted.

We see from the existence of the radius $r=\frac{r_1+r_2}{2}$ that $|\nabla\rho|$ is bounded in $\mathbb{D}$, and so the proof of the theorem is complete.

\section{Lipschitz continuity for certain $(K, K')$-quasiconformal self-mappings of the unit disk}\label{sec-4}
This section is devoted to the proof of Theorem 1.2 together with the statement and the proof of Example 4.1. We start with a lemma.

\subsection{A lemma} $\,$
 By Theorem 1.1, the following assertions easily follows from Lemmas $2.7$ and $2.8$ in \cite{KP}.

\newtheorem*{lemc}{Lemma C}
\begin{lemc}
Under the assumptions of Theorem 1.2, we have

(i) $\lim_{r\to 1^-}\nabla \mathcal{G}(re^{i\theta})=\nabla \mathcal{G}(e^{i\theta})$;

(ii) $\lim_{r\to 1^-}\nabla \mathcal{P}(re^{i\theta})=\nabla \mathcal{P}(e^{i\theta})$ a.e. in $[0, 2\pi]$;

(iii) $\lim_{r\to 1^-}\nabla w(re^{i\theta})=\nabla w(e^{i\theta})$ a.e. in $[0, 2\pi]$;

(iv) $\max_{z\in \mathbb{D}}\{|\mathcal{G}_z(re^{i\theta})|,\;| \mathcal{G}_{\overline{z}}(re^{i\theta})|\}\leq \frac{1}{3}|g|_{\infty}$;

(v) $\max_{0\leq \theta\leq 2\pi}\{|\mathcal{G}_z(e^{i\theta})|,\;|\mathcal{G}_{\overline{z}}(e^{i\theta})|\}\leq \frac{1}{4}|g|_{\infty}$; and

(vi) $\left|J_w(e^{i\theta})- \frac{ \psi'(\theta)}{2\pi}\int_{0}^{2\pi}\frac{|f(e^{i\theta})-f(e^{i\varphi})|^2}{|e^{i\theta}-e^{i\varphi}|^{2}}d\varphi\right|\leq
 \frac{1}{2}\psi'(\theta)|g|_{\infty}$,

where $w|_{\mathbb{S}}=f$ and $f(e^{i\theta})=e^{i \psi(\theta)}$.
\end{lemc}

\subsection{The proof of Theorem 1.2} First, we prove the first statement of the theorem, i.e. the Lipschitz continuity of $w$ in $\mathbb{D}$.
Let $$M_{1}={\rm sup}_{z\in \mathbb{D}}|\nabla w(z)|.$$ Obviously, we only need to show that $M_1$ has an upper bound.
Since $w+\mathcal{G}$ is harmonic, we see from \cite[Lemma 2.2]{Kalaj} and Lemma C that for all $z\in \mathbb{D}$,
\begin{eqnarray*}
|\nabla w(z)| &\leq  &\Big |\nabla w(z)+\nabla \mathcal{G}(z)\Big|+|\nabla \mathcal{G}(z)|\\
&\leq& {\rm esssup}_{0\leq \theta \leq 2\pi}|\nabla w(e^{i \theta})+\nabla \mathcal{G}(e^{i \theta})|+|\nabla \mathcal{G}(z)|\\
&\leq& {\rm esssup}_{0\leq \theta \leq 2\pi}|\nabla w(e^{i \theta})|+\frac{7}{6}|g|_{\infty},
 \end{eqnarray*}
which implies that for every $\varepsilon>0$, there exists a $\theta_{\varepsilon}$ such that
$$
M_{1}
\leq  (1+\varepsilon)|\nabla w(e^{i \theta_{\varepsilon}})| +\frac{7}{6}|g|_{\infty}.
\eqno{(4.1)}$$
Obviously, to estimate $M_1$, it is sufficient to estimate the quantity $|\nabla w(e^{i \theta_{\varepsilon}})|$.

Now, we are going to estimate $|\nabla w(e^{i \theta_{\varepsilon} })|$. It follows from Lemma C and
the assumption $w$ being a $(K,K')$-quasiconformal mapping that
$$|\nabla w(e^{i \theta})|^2=\lim_{r\rightarrow 1^{-}}|\nabla w(re^{i\theta})|^2 \leq K\lim_{r\rightarrow1^{-}}J_{w}(re^{i\theta})+K'\;\;\;\;\;\;\;\;\;\;\;\;\eqno{(4.2)}$$
$$\;\;\;\;\;\;\;\;\;\;\;\;\;\;\;\;\;\;\;\;\;\;\;\;
\leq
 K\psi'(\theta)\Big( \frac{1}{2\pi}\int_{0}^{2\pi}\frac{|f(e^{i\theta}) -f(e^{i \varphi })  |^2}{ |e^{ i\theta }-e^{ i \varphi }|^2 } d \varphi+ \frac{1}{2}|g|_{\infty}\Big) +K'
$$
a.e. in $[0,2\pi]$.

Now, we need a relationship between $\psi'(\theta)$ and $|\nabla w(e^{i\theta})|$. Since Theorem 1.1 guarantees that $|\nabla w(z)|$ is bounded by a constant in $\mathbb{D}$, we deduce from the Lebesgue Dominated
 Convergence Theorem that
\begin{eqnarray*}
f(e^{i \theta})&=&\lim_{r\to 1^-}w(re^{i \theta})=\lim_{r\to 1^-}\int_{\theta_0}^{\theta}\frac{\partial}{\partial \varphi}w(re^{i \varphi})d\varphi+ f(e^{i \theta_{0}})\\
&=&\int_{\theta_0}^{\theta}\lim_{r\to 1^{-}}\Big(\frac{\partial}{\partial \varphi}w(re^{i \varphi})\Big)d\varphi+ f(e^{i \theta_{0}})\\
&=&\int_{\theta_0}^{\theta}\lim_{r\to 1^{-}}\Big(r \nabla w(re^{i \varphi}) ie^{i \varphi}\Big) d\varphi+ f(e^{i \theta_{0}}).
\end{eqnarray*}
Since $f(e^{i\theta})=e^{i\psi(\theta)}$ is absolutely continuous, by differentiating in $\theta$, we obtain
$$
 \frac{d}{d \theta}f(e^{i \theta})=\lim_{r\to 1^{-}}\frac{\partial}{\partial \theta}w(re^{i \theta})=\lim_{r\to 1^{-}}\big(r \nabla w(re^{i \theta})\big) ie^{i \theta}
\eqno{(4.3)}$$
 and
$$
 \frac{d}{d \theta}f(e^{i \theta})=i \psi'(\theta)e^{i \psi(\theta)}\eqno{(4.4)}$$
 a.e. in $[0, 2\pi]$, whence combining Lemma C, we have
$$ \psi'(\theta)= \lim_{r\to 1^{-}}|\nabla w(re^{i \theta})|=|\nabla w(e^{i\theta})|,\eqno{(4.5)}$$ which is our desired relationship between $\psi'(\theta)$ and $|\nabla w(e^{i\theta})|$.

Using (4.5), the relation (4.2) is changed into the following form:
$$
|\nabla w(e^{i \theta})|^2
\leq  K |\nabla w(e^{i\theta})|\Big( \frac{1}{2\pi}\int_{0}^{2\pi}\frac{|f(e^{i \theta }) -f(e^{i\varphi})  |^2}{ |e^{ i\theta }-e^{ i \varphi }|^2 } d\varphi + \frac{1}{2}|g|_{\infty}\Big)+K',
$$
which, necessarily, implies that
$$ |\nabla w(e^{i \theta})|\leq  \frac{K }{2\pi} \int_{0}^{2\pi}\frac{|f(e^{i \varphi}) -f(e^{i \theta})  |^2}{ |e^{ i \varphi }-e^{ i \theta}|^2 }d \varphi+\frac{K }{2 } |g|_{\infty} +\sqrt{K'},  $$
and thus we easily know from
(4.1) that
$$
M_{1}\leq
(1+\varepsilon)\Big[\frac{K}{2\pi}\int_{0}^{2 \pi} M_{1}^{1-\mu}  |e^{i\theta_{\varepsilon}} -e^{i\varphi}|^{ \mu^2+\mu-2} \frac{|f(e^{i\theta_{\varepsilon}})-f(e^{i\varphi})|^{1+\mu}}{|e^{i\theta_{\varepsilon}}- e^{i\varphi}|^{\mu^2+\mu}} \cdot\eqno{(4.6)}$$
$$\;\;\;\;\;\;\;\;\;\;\;\;
\frac{|f(e^{i\theta_{\varepsilon}})-f(e^{i\varphi})|^{ 1-\mu}}{M_{1}^{1-\mu}} d\varphi\Big]
+(1+\varepsilon)\Big[\frac{K}{2}|g|_{\infty}+\sqrt{K'}\Big]+\frac{7}{6}|g|_{\infty},
$$
where $$\mu=\frac{1}{K(1+\pi)^2}.\eqno{(4.7)}$$

Now, we need an auxiliary result which is (4.8) below:
Since for any $z_{1}=re^{i \theta}$ and $z_2=re^{i \eta}$ in $\mathbb{D}$,
 $$|w(z_1)-w(z_2)|\leq \int_{[z_1,z_2]}|\nabla w(z)||dz|\leq M_{1}|z_1 -z_2|,$$
by letting $r\to 1^-$, we obtain
$$ |f(e^{i \theta})-f(e^{i \eta})|\leq  M_{1}| e^{i \theta} - e^{i \eta} |.\eqno{(4.8)}$$

Let us continue the proof. Since $\mathbb{S}$ enjoys the $\frac{\pi}{2}$-chord-arc condition,  we get from (4.6) along with (4.8) and \cite[Lemma 2.4]{KM} that
$$
M_{1}
  \leq
(1+\varepsilon)\Big[\frac{1}{2\pi}KP_{\mathbb{S}}^{1+\mu}\int_{0}^{2 \pi} M_{1}^{1-\mu} |e^{i\theta_{\varepsilon}} -e^{i\varphi}|^{ \mu^2-1}
d\varphi\Big]+(1+\varepsilon)\Big[\frac{1}{2}K|g|_{\infty}+\sqrt{K'}\Big]+\frac{7}{6}|g|_{\infty},
$$
where
$$ P_{\mathbb{S}}=4(1+\pi)2^{\mu}\sqrt{\max\Big\{\frac{2\pi^2 K}{\log 2},\frac{2\pi K'}{K(1+\pi)^2+4}\Big\}} .\eqno{(4.9)}$$
Hence $$
M_{1}
\leq (1+\varepsilon)\Big[M_2M_{1}^{1-\mu}\Big]+(1+\varepsilon)\Big[\frac{1}{2}K|g|_{\infty}+\sqrt{K'}\Big]+\frac{7}{6}|g|_{\infty},
$$
where $$M_2= \frac{1}{2\pi}K P_{\mathbb{S}}^{1+\mu} \int_{0}^{2 \pi}  |e^{i\theta_{\varepsilon}} -e^{i\varphi}|^{ \mu^2-1}  d\varphi.$$
For the convergence of the integral $\int_{0}^{2 \pi}  |e^{i\theta_{\varepsilon}} -e^{i\varphi}|^{ \mu^2-1}  d\varphi$, the reader is referred to \cite[Lemma 1.6]{K2}.
Also, we easily know that $M_2$ does not depend on $\theta_{\epsilon}$. By letting $\varepsilon\to 0$, we get
$$ M_{1}\leq M_2M_{1}^{1-\mu} +\Big(\frac{1}{2}K+\frac{7}{6}\Big)|g|_{\infty}+\sqrt{K'}.\eqno{(4.10)}$$

To get an estimate on $M_1$, we need a lower bound on $M_1$ which is (4.11) below. Since $$\int_{0}^{2\pi}\psi'(\theta)d \theta=\psi (2\pi)-\psi (0)=2\pi,$$
we know that $${\rm esssup}_{0\leq \theta \leq  2\pi}\psi'(\theta)\geq 1.$$

Since $$\psi'(\theta)=|\psi'(\theta)|=\left|\frac{\partial f(e^{i\theta})}{\partial e^{i\theta}}\frac{\partial e^{i\theta}}{\partial\theta}\right|=\left|\frac{\partial f(e^{i\theta})}{\partial e^{i\theta}}\right|$$
and
$${\rm esssup}_{0\leq \theta \leq  2\pi}\lim_{\eta\rightarrow\theta}\left|\frac{f(e^{i\eta})-f(e^{i\theta})}{e^{i\eta}- e^{i\theta}}\right|\leq{\rm esssup}_{0\leq \theta \not= \eta< 2\pi}\left|\frac{f(e^{i\eta})-f(e^{i\theta})}{e^{i\eta}- e^{i\theta}}\right|,$$
it follows from the inequality (4.8) that
$$1\leq {\rm esssup}_{0\leq \theta \leq  2\pi}\psi'(\theta)\leq {\rm esssup}_{0\leq \theta \not= \eta< 2\pi}\left|\frac{f(e^{i\eta})-f(e^{i\theta})}{e^{i\eta}- e^{i\theta}}\right|\leq M_{1}.\eqno{(4.11)}$$

Now, we are able to get an upper bound for $M_1$. Using (4.11), the relation (4.10) implies
$$ M_{1}\leq \Big(M_2 +\Big(\frac{1}{2}K+\frac{7}{6}\Big)|g|_{\infty}+\sqrt{K'}\Big)M_{1}^{1-\mu},$$
and so
$$
M_{1} \leq \Big( M_2  +\frac{1}{2}K|g|_{\infty}+\frac{7}{6} |g|_{\infty}+\sqrt{K'} \Big)^{K(1+\pi)^2}
= C_0.
$$

Moreover, by \cite[Lemma 2.9]{KP} and (4.10), we see that if
$$(1-\mu) M_2=\frac{1}{2\pi}\Big(1-\frac{1}{K(1+\pi)^2}\Big)K P_{\mathbb{S}}^{1+\mu} \int_{0}^{2 \pi}  |e^{i\theta_{\varepsilon}} -e^{i\varphi}|^{ \mu^2-1}  d\varphi
<1,
$$
then $$M_{1}\leq \frac{M_2+\frac{1}{2}K|g|_{\infty}+\frac{7}{6} |g|_{\infty}+\sqrt{K'}-(1-\mu) M_2}{1-(1-\mu) M_2}=C_1.$$
Let
$$
M=\begin{cases}
\displaystyle C_0, \;\;\;\;\;\;\;\;\;\;\;\;\;\;\;\;\, {\rm if} \;\, (1-\mu) M_2\geq 1, \\
\displaystyle \min\{C_0, C_1\},\,\; {\rm if} \;\, (1-\mu) M_2< 1.
\end{cases}$$ Then we see that
$$M_1\leq M,$$ and so the proof of the first statement of the theorem is complete.
\medskip

Next, we are going to prove the second statement of the theorem, i.e. the coLipschitz continuity of $w$ under the assumption that $w^{-1}$ is also $(K, K')$-quasiconformal.
It follows from Lemma C that
$$l(\nabla w)\geq l(\nabla \mathcal{P})-|\nabla \mathcal{G}|\geq  l(\nabla \mathcal{P})- \frac{2}{3}|g|_{\infty}\eqno{(4.12)}$$ a.e. in $\mathbb{D}$.

Obviously, to prove the coLipschitz continuity of $w$, it is sufficient to find the lower bound of $l(\nabla w)$ in $\mathbb{D}$, and (4.12) implies that it is enough to find the lower bound of $l(\nabla \mathcal{P})$.
For this purpose, we need the following claim.

\vspace{0.3truecm}
\noindent{\bf Claim 4.1.}
For a.e. $\theta\in [0,2\pi]$, we have
$$K\psi'(\theta)\geq \frac{1}{2\pi}\int_{0}^{2\pi}\frac{|f(e^{i\theta})-f(e^{i\varphi})|^2}{|e^{i\theta}-e^{i\varphi}|^2}d\varphi
-\frac{1}{2}|g|_{\infty}-\sqrt{K'}\geq N_1,$$
 where
 $$N_1= \frac{1}{2\pi} P^{-\frac{2}{\mu}}_{\mathbb{S}}\int_{0}^{2\pi}|e^{i\theta}-e^{i\varphi}|^{\frac{2}{\mu}-2}d\varphi-\frac{1}{2}|g|_{\infty}-\sqrt{K'} .$$

First, we prove the inequality:
$$ K\psi'(\theta)\geq \frac{1}{2\pi}\int_{0}^{2\pi}\frac{|f(e^{i\theta})-f(e^{i\varphi})|^2}{|e^{i\theta}-e^{i\varphi}|^2}d\varphi
-\frac{1}{2}|g|_{\infty}-\sqrt{K'}\eqno{(4.13)}$$  a.e. in $[0,2\pi]$.

It follows from
(4.3) and (4.4) that $$\psi'(\theta)=\left|\frac{d f(e^{i\theta})}{d \theta} \right|=\left| \lim_{r\to1^{- }} \frac{\partial w(re^{i\theta})}{\partial \theta} \right|.$$
Since  Lemma C guarantees that
$$\psi'(\theta)=\left| \lim_{r\to1^{-}} \frac{\partial w(re^{i\theta})}{\partial \theta} \right|\geq |w_{z}(e^{i\theta})|-|w_{\overline{z}}(e^{i\theta})|=l(\nabla w(e^{i\theta}))$$
and
$$\frac{J_{w}(e^{i\theta})}{\psi'(\theta)}\geq \frac{1}{2\pi}\int_{0}^{2\pi}\frac{|f(e^{i\theta})-f(e^{i\varphi})|^2}{|e^{i\theta}-e^{i\varphi}|^2}d\varphi-\frac{1}{2}|g|_{\infty}$$ a.e. in $[0,2\pi]$,
we have
\begin{eqnarray*}
\frac{1}{2\pi}\int_{0}^{2\pi}\frac{|f(e^{i\theta})-f(e^{i\varphi})|^2}{|e^{i\theta}-e^{i\varphi}|^2}d\varphi-\frac{1}{2}|g|_{\infty}
  &\leq &\frac{J_{w}(e^{i\theta})}{\psi'(\theta)}\leq |\nabla w(e^{i\theta})|\\ \nonumber
  &\leq & K l(\nabla w(e^{i\theta}))+\sqrt{K'}\\ \nonumber
  &\leq & K \psi'(\theta)+\sqrt{K'}
\end{eqnarray*}  a.e. in $[0,2\pi]$, as required.\medskip

Next, we get an estimate on the integral in Claim 4.1, which is as follows:
$$
\frac{1}{2\pi}\int_{0}^{2\pi}\frac{|f(e^{i\theta})-f(e^{i\varphi})|^2}{|e^{i\theta}-e^{i\varphi}|^2}d\varphi
\geq  \frac{1}{2\pi} P^{-\frac{2}{\mu}}_{\mathbb{S}}\int_{0}^{2\pi}|e^{i\theta}-e^{i\varphi}|^{\frac{2}{\mu}-2}d\varphi.
\eqno{(4.14)}$$

By \cite[Lemma 2.4]{KM} and the assumptions that $w^{-1}$ is $(K,K')$-quasiconformal and $w$ is normalized, we have that for all $z_1$ and $z_2\in \mathbb{S}$,
$$|z_1-z_2|\leq P_{\mathbb{S}}|w (z_1)-w (z_2)|^{\mu},$$
i.e. $$|w(z_1)-w(z_2)|\geq P^{-\frac{1}{\mu}}_{\mathbb{S}} | z_1 - z_ 2|^{\frac{1}{\mu}},$$
where $P_{\mathbb{S}}$ and $\mu$ are the same as in (4.7) and (4.9).
Then we have
$$
\frac{1}{2\pi}\int_{0}^{2\pi}\frac{|f(e^{i\theta})-f(e^{i\varphi})|^2}{|e^{i\theta}-e^{i\varphi}|^2}d\varphi
 \geq N_2,
$$
where $$N_2=\frac{1}{2\pi} P^{-\frac{2}{\mu}}_{\mathbb{S}}\int_{0}^{2\pi}|e^{i\theta}-e^{i\varphi}|^{\frac{2}{\mu}-2}d\varphi
\leq P^{-\frac{2}{\mu}}_{\mathbb{S}}2^{\frac{2}{\mu}-2}.$$
Here we remark that by using the substitution in the integral, we easily see that $N_2$ is independent of $\theta$, i.e. $N_2=N_2(K, K')$.

Obviously, the proof of Claim 4.1 follows from (4.13) and (4.14).
\medskip

Now, we are ready to finish the proof of our theorem by applying Claim 4.1.
It follows from Claim 4.1, together with the inequalities (4.5), Lemmas 1.1 and C, that
$$ N_1\leq K\psi'(\theta) = K|\nabla w(e^{i\theta})| \leq K^2l(\nabla w(e^{i\theta}))+K\sqrt{K'}$$
A.E. In $[0, 2\Pi]$, Whence Again Lemma C Implies
$$
l(\nabla \mathcal{P}(e^{i\theta}))=\lim_{r\to 1^-}l(\nabla \mathcal{P}(re^{i\theta}))=\lim_{r\rightarrow 1^-} (|\mathcal{P}_z(re^{i\theta})|-|\mathcal{P}_{\overline{z}}(re^{i\theta})|)\;\;\;\;\;\;\;\;\;\;
\eqno{(4.15)}$$
\begin{eqnarray*}
&\geq& \lim_{r\rightarrow1^-} (|w_z(re^{i\theta})+\mathcal{g}_z(re^{i\theta})|-|w_{\overline{z}}(re^{i\theta})+\mathcal{g}_{\overline{z}}(re^{i\theta})|)
\\
&\geq & \lim_{r\rightarrow 1^-}\Big(l(\nabla w(re^{i\theta}))-|\nabla \mathcal{g}(re^{i\theta})|\Big)\\
&\geq & \frac{N_1}{K^2 }-\frac{\sqrt{K'}}{K}-\frac{1}{2}|G|_{\infty}.
\end{eqnarray*}

Let $$N=\frac{N_1}{K^2}-\frac{\sqrt{K'}}{K}-\frac{7}{6}|G|_{\infty}.$$ Then $N=N(K, K', |G|_{\infty})$ Since $N_1=N_1(K, K', |G|_{\infty})$, And Next, We Are Going To Show The Following.

\vspace{0.3truecm}
\noindent{\bf Claim 4.2.}
The inequality $l(\nabla w) \geq  N$ holds in $\mathbb{D}$.

Without loss of generality, we assume that $N>0$. Obviously, $$\frac{N_1}{K^2}-\frac{\sqrt{K'}}{K}-\frac{1}{2}|g|_{\infty}>0.$$

Under this assumption, we need to get a lower bound for $l(\nabla \mathcal{P})$ in $\mathbb{D}$ (See (4.16) below). We will employ the famous Heinz Theorem \cite{he} to reach this aim.

Since $f(e^{i\theta})=e^{i \psi(\theta)}$ is an increasing homeomorphism on $ \mathbb{S}$, we see from the Choquet-Rad\'{o}-Kneser Theorem (cf. \cite[p. 29]{D}) that $\mathcal{P}$ is a sense-preserving harmonic diffeomorphism. Then Heinz Theorem implies
$$2|\mathcal{P}_z|^2\geq|\mathcal{P}_z|^2+|\mathcal{P}_{\overline{z}}|^2\geq \frac{1}{\pi^2}.$$
Let $$\varphi(z)= \frac{\overline{\mathcal{P}_{\overline{z}}}(z)}{\mathcal{P}_z(z)}\;\;\text{and}\;\;\phi(z)=\frac{1}{\mathcal{P}_z(z)} \Big(  \frac{N_1}{K^2}-\frac{\sqrt{K'}}{K}-\frac{1}{2}|g|_{\infty}\Big).$$  Then both $\varphi$ and $\phi$ are holomorphic, $|\varphi(z)|<1$ and further
$$|\phi|\leq \sqrt{2}\pi \Big(  \frac{N_1}{K^2}-\frac{\sqrt{K'}}{K}-\frac{1}{2}|g|_{\infty}\Big)$$ in $\mathbb{D}$.
Since (4.15) leads to
$$|\varphi(e^{i\theta})|+|\phi(e^{i\theta})|
=\frac{|\mathcal{P}_{\overline{z}}(e^{i\theta})|+\big(\frac{N_1}{K^2}-\frac{\sqrt{K'}}{K}-\frac{1}{2}|g|_{\infty}\big)}{|\mathcal{P}_z(e^{i\theta})|}
\leq 1,$$
we see that
$$
|\varphi(z)|+|\phi(z)|\leq P[|\varphi|_{\mathbb{S}}|](z)+P[|\phi|_{\mathbb{S}}|](z)\leq 1,
$$
which implies
$$|\mathcal{P}_{\overline{z}}(z)|+\frac{N_1}{K^2}- \frac{\sqrt{K'}}{K}-\frac{1}{2}|g|_{\infty}\leq|\mathcal{P}_z(z)|$$ in $ \mathbb{D}$,
i.e. $$l(\nabla \mathcal{P})\geq\frac{N_1}{K^2}- \frac{\sqrt{K'}}{K}-\frac{1}{2}|g|_{\infty}\eqno{(4.16)}$$ in $\mathbb{D}$, as required.\medskip

Now, it follows from (4.12) and (4.16) that
$$l(\nabla w) \geq  N,$$
and hence the claim is proved. \medskip

Since the second statement in Theorem 1.2 $(2)$ easily follows from Claim 4.2, we see that the proof of the theorem is complete.
\qed

\subsection{An example} The following example shows that the assumption ``$w^{-1}$ being $(K,K')$-quasiconformal" in the second assertion in Theorem 1.2 is necessary.

\newtheorem*{exa4.1}{Example 4.1}
\begin{exa4.1} Let $w(z)=\frac{1}{2n}\big((2n+1)z-z|z|^{2n}\big)$ in $\mathbb{D}$, where $n\geq 1$. Then

(1) $w$ is a $(K,K')$-quasiconformal mapping of $\mathbb{D}$ onto $\mathbb{D}$;

(2) $w$ is not $K$-quasiconformal for any $K\geq1$;

(3) $w^{-1}$ is not a $(K,K')$-quasiconformal mapping for any $K\geq 1$ and $K'\geq 0$;

(4) $w$ is Lipschitz continuous but not coLipschitz.

\end{exa4.1}
\begin{proof}
Elementary calculations show that
$$w_z(z)=\frac{1}{2n}\big((2n+1)-(n+1)|z|^{2n}\big)\;\;\text{and}\;\;w_{\overline{z}}(z)= -\frac{1}{2}|z|^{2n-2}z^2.$$
It follows from $$J_{w}(z)=\frac{1}{4n^2}\big(|w_z(z)|^2-|w_{\overline{z}}(z)|^2\big)=\frac{1}{4n^2} (2n+1)(1-|z|^{2n})(2n+1-|z|^{2n}) >0,$$
together with the fact $w(e^{i \theta})  =e^{i \theta}$ and the degree principle, that $w$ is a sense-preserving homeomorphism from $\mathbb{D}$ onto $\mathbb{D}$. Furthermore,
 $$l( \nabla w)(z)=\frac{1}{2n}\big((2n+1)-(2n+1)|z|^{2n}\big)\;\;\text{and}\;\;|\nabla w(z)|=\frac{1}{2n}\big((2n+1)-|z|^{2n}\big).$$ Then $w$ is a $\big(1,(1+\frac{1}{2n})^2\big)$-quasiconformal mapping since
 $$|\nabla w |^2\leq J_w +|\nabla w |^2\leq J_{w}+\big(1+\frac{1}{2n}\big)^2.$$ The limit
 $$\lim_{|z|\to 1^{-}} \frac{|\nabla w |^2}{J_{w}}=+\infty$$ tells us that
  $w$ is not a $K$-quasiconformal mapping for any $K\geq 1$. The Lipschitz continuity of $w$ easily follows from the estimate $|\nabla w|\leq 1+\frac{1}{2n}$.

   It is well known that for the nonsingular matrix $\nabla w$, we have
   $$\big|\nabla  w^{-1}\big|=1/l( \nabla w)\;\; \mbox{and}\;\; l(\nabla w^{-1})=1/|\nabla w|\;\, {\rm (cf.}\; {\rm [11])}.$$ Hence
   $$|\nabla w^{-1}(z) | =\frac{2n}{(2n+1)-(2n+1)|z|^{2n}}\;\;\text{and}\;\;l(\nabla w^{-1}(z)) =\frac{2n}{(2n+1)- |z|^{2n}}.$$
  Then for any $K\geq 1$, $$\lim_{|z|\to 1^{-}} \Big( \big|\nabla w^{-1} (z) \big|^2-K J_{w^{-1}}(z)\Big) =\lim_{|z|\to 1^{-}}
  |\nabla  w^{-1}|\big(|\nabla  w^{-1}|-K  l( \nabla w^{-1})\big)
  =+\infty.$$
 This shows that $w^{-1}$ is not a $(K,K')$-quasiconformal mapping.

  Let $\partial_{\alpha}w(z)$ denote the directional derivative of $w$.
  Note that if $w$ is Lipschitz continuous with Lipschitz constant $C$, then
  $$\big | \partial_{\alpha}w(z) \big|=\left|\lim_{r\to 0}\frac{w(z+re^{i \alpha})-w(z)}{r } \right|= \lim_{r\to0}\frac{\left|w(z+re^{i \alpha})-w(z)\right|}{r } \leq C.\eqno{(4.17)}$$ Hence it follows from the obvious fact $|\nabla w(z)|=\max_{\alpha}\big | \partial_{\alpha}w(z) \big|$ that $$|\nabla w|\leq C.$$ That's, ``$w$ being Lipschitz continuous" is equivalent to ``$|\nabla w|$ being bounded".
  Since $$\lim_{|z|\to1^{-}}\big|\nabla (w^{-1} )(z)\big| =+\infty,$$
we see that $w^{-1}$ is not Lipschitz continuous and so $w $ is not coLipschitz continuous. The proof is finished.
\end{proof}

\end{document}